\documentclass[reqno,a4paper,11pt]{amsart}

\usepackage{amsmath,amstext,amsfonts,amsbsy,eucal,amssymb}
\usepackage[latin1]{inputenc}

\numberwithin{equation}{section}

\newtheorem{theorem}{Theorem}[section]
\newtheorem{lemma}[theorem]{Lemma}

\newtheorem{proposition}[theorem]{Proposition}

\theoremstyle{definition}

\begin{document}

\parskip 4pt
\baselineskip 16pt


\title[Bernoulli numbers and sums of powers of integers of higher order]
{Bernoulli numbers and sum of powers of integers of higher order}

\author[Andrei K. Svinin]{Andrei K. Svinin}
\address{Andrei K. Svinin, 
Matrosov Institute for System Dynamics and Control Theory of 
Siberian Branch of Russian Academy of Sciences,
P.O. Box 292, 664033 Irkutsk, Russia}
\email{svinin@icc.ru}

\author[Svetlana V. Svinina]{Svetlana V. Svinina}
\address{Svetlana V. Svinina, 
Matrosov Institute for System Dynamics and Control Theory of 
Siberian Branch of Russian Academy of Sciences,
P.O. Box 292, 664033 Irkutsk, Russia}
\email{gaidamak@icc.ru}

%
%

\date{\today}



\begin{abstract}
We give an expression of polynomials for higher sums of powers of integers via the higher order Bernoulli numbers.
\end{abstract}

\maketitle

\section{Introduction}

As is known the sum of powers of integers  \cite{Graham}, \cite{Knut}
\[
S_m(n):=\sum_{q=1}^{n}q^{m}
\]
can be computed with the help of some appropriate polynomial $\hat{S}_m(n)$ for any $m\geq 0$. Exponential generating function for  the sums $S_m(n)$ is given by
\begin{equation}
S(n, t)=\sum_{q=1}^{n}e^{qt}=\frac{e^{(n+1)t}-e^t}{e^t-1}.
\label{genf}
\end{equation}
Expanding in series (\ref{genf}) yields an infinite set of polynomials $\{\hat{S}_m(n) : m\geq 0\}$, that is,
\[
S(n, t)=\sum_{q\geq 0}\hat{S}_{q}(n)\frac{t^q}{q!}.
\]
It is a classical result  that these polynomials can be expressed as \cite{Jacobi}
\begin{equation}
\hat{S}_{m}(n)=\frac{1}{m+1}\sum_{q=0}^m(-1)^q{m+1\choose q}B_{q}n^{m+1-q},
\label{2}
\end{equation}
where $B_q$ are the Bernoulli numbers that can be derived from the exponential  generating function
\begin{equation}
\frac{t}{e^t-1}=\sum_{q\geq 0}B_q\frac{t^q}{q!}.
\label{Bernoulli}
\end{equation}
It follows from (\ref{Bernoulli})  that the Bernoulli numbers satisfy  the recurrence relation
\begin{equation}
\sum_{q=0}^{m}{m+1\choose q}B_{q}=\delta_{0, m}.
\label{rec-rel}
\end{equation}
This relation is in fact  the simplest one of many known  recurrence relations involving the Bernoulli numbers (see, for example, \cite{Agoh} and references therein). One can derive, for example, an  infinite number of  recurrence relations of the form
\begin{equation}
\sum_{q=0}^{m}{m+k\choose q}S(m+k-q, k)B_{q}=\frac{m+k}{k}S(m+k-1, k-1),\;\;
\forall k\geq 1.
\label{rec-rel1}
\end{equation}

In this paper we investigate a class of sums that correspond to a $k$-th power of generating function (\ref{genf}) for $k\geq 1$. Our main result is a formula for  polynomials allowing to calculate   these sums. It turns out that these polynomials are expressed via the higher Bernoulli numbers.

\section{The power sums of higher order}

Let us now consider a power of the generating function (\ref{genf}):
\[
\left(S(n, t)\right)^k:=\sum_{q\geq 0}S_{q}^{(k)}(n)\frac{t^q}{q!}.
\]
We have
\begin{equation}
\left(S(n, t)\right)^k=\left(\sum_{q=1}^{n}e^{qt}\right)^k=\sum_{q=k}^{kn}{k\choose q}_ne^{qt}.
\label{1}
\end{equation}
The coefficients ${k\choose q}_n$  obviously generalizing the binomial coefficients  originated from Abraham De Moivre and Leonard Euler works \cite{Moivre}, \cite{Euler} and extensively studied in the literature due to their  applicability. From (\ref{1}), we see that it is a generating function for the sums of the form
\begin{equation}
S_{m}^{(k)}(n):=\sum_{q=0}^{k(n-1)}{k\choose q}_{n}\left(k+q\right)^m.
\label{sums} 
\end{equation}
It is natural to call (\ref{sums}) the  sums of powers of integers of higher order. Expanding
\[
\left(\frac{e^{(n+1)t}-e^t}{e^t-1}\right)^k=\sum_{q\geq 0}\hat{S}_{q}^{(k)}(n)\frac{t^q}{q!},
\]
we get an infinite number of polynomials $\hat{S}_{m}^{(k)}(n)$. 

Our goal in the paper is to prove that 
\begin{equation}
\hat{S}_{m}^{(k)}(n)=\frac{1}{{m+k\choose k}}\sum_{q=0}^{m}(-1)^q{m+k\choose q}B_q^{(k)}S(m+k-q, k)n^{m+k-q},
\label{higher-polynomials}
\end{equation}
where $B_q^{(k)}$ are the higher order Bernoulli numbers  defined as 
\begin{equation}
\frac{t^k}{(e^t-1)^k}=\sum_{q\geq 0}B^{(k)}_q\frac{t^q}{q!}.
\label{Bern-high}
\end{equation}
The Bernoulli numbers of higher order appeared  in \cite{Norlund} in connection with a theory of finite differences   
and then was investigated by many authors from different points of view (see, for example, \cite{Carlitz}). These numbers are known to satisfy \cite{Norlund}
\[
B_n^{(k+1)}=\frac{k-n}{k}B_n^{(k)}-nB_{n-1}^{(k)}.
\]
The number $B^{(k)}_n$ with fixed $n\geq 0$ turns out to be some polynomial in $k$. These kind of polynomials are known as  N\"orlund polynomials.  One can find a number of these polynomials in  \cite{Norlund}. For convenience, we have written out several N\"orlund  polynomials in the Appendix. The numbers $S(n, k)$ in (\ref{higher-polynomials}) are the Stirling numbers of the second kind that satisfy recurrence relation 
\begin{equation}
S(n, k)=S(n-1, k-1)+kS(n-1, k)
\label{rr}
\end{equation}
with appropriate boundary conditions \cite{Weisstein2}, \cite{Graham}.

It is easy to prove that the higher order Bernoulli numbers satisfy the recurrence relation
\begin{equation}
\sum_{q=0}^{m}{m+k\choose q}S(m+k-q, k)B_{q}^{(k)}=\delta_{0, m}.
\label{impl}
\end{equation}
The most general relation involving (\ref{rec-rel}), (\ref{rec-rel1}) and (\ref{impl}) as particular cases is
\begin{equation}
\sum_{q=0}^{m}{m+k\choose q}S(m+k-q, k)B_{q}^{(r)}=\frac{{m+k\choose k}}{{m+k-r\choose k-r}}S(m+k-r, k-r),\;\; \forall k\geq r.
\label{impl1}
\end{equation}

As is known $S(m+k, k)$, for any fixed $m\geq 0$, is expressed as a polynomial  $f_m(k)$ of degree $2m$, which satisfy the identity 
\[
f_m(k)-f_m(k-1)=kf_{m-1}(k)
\]
following  from the identity (\ref{rr}). Therefore we can replace $S(m+k-q, k)$ by $f_{m-q}(k)$ in (\ref{higher-polynomials}). In the literature the polynomials $f_m(k)$ are known as the Stirling polynomials \cite{Gessel}, \cite{Jordan}. These  are known to be expressed via the N\"orlund polynomials as (see, for example, \cite{Adelberg})
\[
f_m(k)={m+k\choose m}B_m^{(-k)}.
\]

The following proposition also gives the relationship of the higher Bernoulli numbers with the Stirling numbers.
\begin{proposition}
One has
\begin{equation}
B_m^{(k)}=\sum_{q=1}^{m} \frac{s(q+k, k)}{{q+k\choose k}}S(m, q).
\label{8}
\end{equation}
\end{proposition}
In (\ref{8}), $s(n, k)$ stands for the Stirling numbers of the first kind \cite{Weisstein1}.

It is evident that in the case $k=1$, (\ref{8}) becomes 
\[
B_m=\sum_{q=1}^{m}(-1)^{q} \frac{q!}{q+1}S(m, q),
\]
while  in the case $k=2$, it takes the following form:
\[
B_m^{(2)}=2\sum_{q=1}^{m}(-1)^{q} \frac{(q+1)!H_{q+1}}{(q+1)(q+2)}S(m, q),
\]
where $H_m$ are harmonic number defined by $H_m:=\sum_{q=1}^{m}1/q$.

To prove (\ref{higher-polynomials}), we need the following lemma: 
\begin{lemma} \label{le2}
By virtue of (\ref{impl1}) we have
\begin{eqnarray}
R_m^{(k, r)}(n)&:=&\sum_{q=0}^{m}(-1)^{q}{m+k\choose q}S(m+k-q, k)\hat{S}_q^{(r)}(n)\label{le1}\\
      &=&\frac{1}{{k\choose r}}\sum_{j=0}^{m}(-1)^j{m+k\choose m+k-r-j}S(m+k-r-j, k-r)\nonumber\\
			&&\times S(r+j, r)n^{r+j},\;\; \forall m\geq 0,\;\;k\geq r.
			\label{lee}
\end{eqnarray}
\end{lemma}
It should be remarked that in the case $k=r$, (\ref{lee}) becomes
\begin{equation}
R_m^{(k, k)}(n)=(-1)^mS(m+k, k)n^{m+k}.
\label{rec-rel2}
\end{equation}

\noindent
\textbf{Proof of lemma \ref{le2}}. We can rewrite (\ref{le1}) as  
\[
R_{m}^{(k, r)}(n)=\sum_{0\leq j\leq q \leq m}a_qb_{q, j}n^{r+q-j},
\]
where
\[
a_{q}:=\frac{{m+k\choose q}}{{r+q\choose r}}S(m+k-q, k)
\]
and
\[
b_{q, j}:=(-1)^{q-j}{r+q\choose j}B_j^{(r)}S(r+q-j, r).
\]
Let $\tilde{j}=q-j$ and 
\[
b_{q, \tilde{j}}=(-1)^{\tilde{j}}{r+q\choose q- \tilde{j}}B_{q- \tilde{j}}^{(r)}S(r+\tilde{j}, r).
\]
In what follows, for simplicity,   let us write $\tilde{j}$ without the tilde. Making use the identity
\[
{r+q\choose q- j}={r+q\choose r+j}={q\choose j}\frac{{r+q\choose r}}{{r+j\choose r}},
\]
we get 
\[
a_{q}b_{q, j}=(-1)^{j}\frac{S(r+j, r)}{{r+j\choose r}}{m+k\choose q}{q\choose j}S(m+k-q, k)B_{q- j}^{(r)}
\]
and therefore 
\begin{eqnarray}
R_{m}^{(k, r)}(n)&=&\sum_{0\leq j\leq q \leq m}a_qb_{q, j}n^{r+j}\nonumber\\ 
        &=&\sum_{0\leq j\leq m}(-1)^{j}\frac{S(r+j, r)}{{r+j\choose r}}n^{r+j}\sum_{j\leq q\leq m}{m+k\choose q}{q\choose j}S(m+k-q, k)B_{q- j}^{(r)}.\nonumber
\end{eqnarray}
In turn, making use  the identity
\[
{m+k\choose q}{q\choose j}={m+k\choose j}{m+k-j\choose q-j},
\]
we get
\begin{eqnarray}
R_{m}^{(k, r)}(n)&=&\sum_{0\leq j\leq m}(-1)^{j}\frac{S(r+j, r)}{{r+j\choose k}}{m+k\choose j}n^{r+j}\nonumber\\
&&\times\sum_{j\leq q\leq m}{m+k-j\choose q-j}S(m+k-q, k)B_{q- j}^{(r)}. \nonumber
\end{eqnarray}
Finally, by virtue of (\ref{impl1}), we get
\begin{eqnarray}
&&\sum_{j\leq q\leq m}{m+k-j\choose q-j}S(m+k-q, k)B_{q- j}^{(r)}\nonumber\\
&&\;\;\;\;\;\;\;\;\;\;\; =\sum_{0\leq q\leq m-j}{m+k-j\choose q}S(m+k-j-q, k)B_{q}^{(r)}\nonumber\\
&&\;\;\;\;\;\;\;\;\;\;\; =\frac{{m+k-j\choose k}}{{m+k-j-r\choose k-r}}S(m+k-j-r, k-r)\nonumber
\end{eqnarray}
and hence
\begin{eqnarray}
R_{m}^{(k, r)}(n)&=&\sum_{0\leq j\leq m}(-1)^{j}\frac{{m+k\choose j}}{{r+j\choose r}}\frac{{m+k-j\choose k}}{{m+k-j-r\choose k-r}}S(m+k-j-r, k-r)\nonumber\\
                 &&\times S(r+j, r)n^{r+j}\nonumber\\
								 &=&\frac{1}{{k\choose r}}\sum_{0\leq j\leq m}(-1)^{j}{m+k\choose m+k-j-r}S(m+k-j-r, k-r).\nonumber\\
								&&\times S(r+j, r)n^{r+j}\nonumber
\end{eqnarray}
Therefore the lemma is proved. $\Box$

The   recurrence relation, for example (\ref{rec-rel2}), uniquely determines an infinite set of polynomials $\{\hat{S}_m^{(k)}(n) : m\geq 0\}$. We have written out some of them in the Appendix. For example,
$\hat{S}_{0}^{(k)}(n)=n^k$. On the other hand
\[
S_{0}^{(k)}(n):=\sum_{q=0}^{k(n-1)}{k\choose q}_{n}=\Biggl(\sum_{q=0}^{n-1}t^q\Biggr)^{k}\Bigl|_{t=1}=n^k.
\]
\begin{lemma} \label{le3}
The higher sums $S_m^{(k)}(n)$ satisfy the same recurrence relations as in lemma \ref{le2}, that is,
\begin{eqnarray}
&&\sum_{q=0}^{m}(-1)^{q}{m+k\choose q}S(m+k-q, k)S_q^{(r)}(n) \nonumber  \\
      &&\;\;\;=\frac{1}{{k\choose r}}\sum_{j=0}^{m}(-1)^j{m+k\choose m+k-r-j}S(m+k-r-j, k-r)S(r+j, r)n^{r+j}.						
			\label{id}
\end{eqnarray}

\end{lemma}
In the case $k=r=1$, (\ref{id}) becomes the well-known identity for the sums of powers \cite{Riordan}.

\noindent
\textbf{Proof of lemma \ref{le3}}. This lemma is proved by using standard arguments. Let us replace an argument  of generating function $t\rightarrow -t$ to get
\begin{equation}
\sum_{q\geq 0}(-1)^qS^{(r)}_{q}(n)\frac{t^q}{q!}=(-1)^k\left(\frac{e^{-nt}-1}{e^{t}-1}\right)^{r}.
\label{id1}
\end{equation}
Multiplying both sides of (\ref{id1}) by $(e^t-1)^k$ and taking into account that
\[
(e^t-1)^k=k!\left(\sum_{q\geq 0} S(q, k)\frac{t^q}{q!}\right),
\]
we get (\ref{id}). $\Box$

Now, we are in a position to prove our theorem.
\begin{theorem}
One has
\[
S_{m}^{(k)}(n)=\hat{S}_{m}^{(k)}(n),
\]
where $\hat{S}_{m}^{(k)}(n)$ be the polynomials (\ref{higher-polynomials}).
\end{theorem}

\noindent
\textbf{Proof}. This theorem is a simple consequence of lemma \ref{le2} and  lemma \ref{le3} since the sums ${S}_m^{(k)}(n)$ satisfy the same recurrence relations  as the polynomials  $\hat{S}_m^{(k)}(n)$. $\Box$

\section{The relationship of the sums $S_{m}^{(k)}(n)$ to other sums}

In \cite{Svinin} we  considered  sums of the form
\begin{equation}
\mathcal{S}_{m}^{(k)}(n):=\sum_{\{\lambda\}\in B_{j, jn}}\left(\lambda_1^{m}+(\lambda_2-n)^{m}+\cdots+(\lambda_k-kn+n)^{m}\right),
\label{sums11} 
\end{equation}
where it is supposed that $m$ is odd. Here $B_{k, kn}:=\{\lambda_q : 1\leq \lambda_1\leq \cdots \leq \lambda_k\leq kn\}$. Let us remark that   there are some terms of the form $r^m$ with negative $r$ in (\ref{sums11}). It is evident that in this case $r^m=-|r|^m$. By this rule, the sum (\ref{sums11}) can be rewritten as
\begin{equation}
\mathcal{S}_{m}^{(k)}(n)=\sum_{q=0}^{kn}c_q(k, n)q^m
\label{defin}
\end{equation}
with some integer coefficients $c_r(k, n)$.

It was conjectured in \cite{Svinin} that in the case of odd $m$ the sums $\mathcal{S}_{m}^{(k)}(n)$ and $S_{m}^{(k)}(n)$ are related with each other by
\begin{equation}
\mathcal{S}_{m}^{(k)}(n)=\sum_{q=0}^{k-1}{k(n+1)\choose q} S_{m}^{(k-q)}(n).
\label{relationsh}
\end{equation}
Let us define the sums $\mathcal{S}_{m}^{(k)}(n)$ with even $m$ by (\ref{defin}). It is evident that conjectural relation (\ref{relationsh}) is valid for both odd  and even $m$. More exactly, actual calculations show that
\[
\mathcal{S}_{m}^{(k)}(n)-\sum_{q=0}^{k-1}{k(n+1)\choose q} S_{m}^{(k-q)}(n)=c_0(k, n)\delta_{m, 0}.
\]
\section*{Appendix}

\subsection*{N\"orlund polynomials}

The first six of the N\"orlund polynomials  are given by
\[
B_0^{(k)}=1,\;\;
B_1^{(k)}=-\frac{1}{2}k,\;\;
B_2^{(k)}=\frac{1}{12}k\left(3k-1\right),\;\;
B_3^{(k)}=-\frac{1}{8}k^2\left(k-1\right),
\]
\[
B_4^{(k)}=\frac{1}{240}\left(15k^3-30k^2+5k+2\right),\;\;
B_5^{(k)}=-\frac{1}{96}k^2\left(k-1\right)\left(3k^2-7k-2\right).
\]

\subsection*{The polynomials $\hat{S}_m^{(k)}(n)$}

The first six of these polynomials are given by 
\[
\hat{S}_0^{(k)}(n)=n^k,\;\;
\hat{S}_1^{(k)}(n)=\frac{k}{2}n^k\left(n+1\right),\;\;
\]
\[
\hat{S}_2^{(k)}(n)=\frac{k}{12}n^k(n+1)\left((3k+1)n+3k-1\right),
\]
\[
\hat{S}_3^{(k)}(n)= \frac{k^2}{8}n^k(n+1)^2\left((k+1)n+k-1\right),
\]
\begin{eqnarray}
\hat{S}_{4}^{(k)}(n)&=&\frac{k}{240}n^k(n+1)\left((15k^3+30k^2+5k-2)n^3+(45k^3+30k^2-5k+2)n^2\right. \nonumber\\
                   &&\left.+(45k^3-30k^2-5k-2)n+15k^3-30k^2+5k+2\right),\nonumber
\end{eqnarray}
\begin{eqnarray}
\hat{S}_{5}^{(k)}(n)&=&\frac{k^2}{96}n^k(n+1)^2\left((3k^3+10k^2+5k-2)n^3+(9k^3+10k^2-5k+2)n^2\right.\nonumber\\
           &&\left.+(9k^3-10k^2-5k-2)n+3k^3-10k^2+5k+2\right).\nonumber
\end{eqnarray}

\end{document}